\documentclass[a4paper,11pt]{article}
\usepackage{amssymb,amsmath, amsthm}
\usepackage{multirow}

\title{Local heights on elliptic curves and intersection multiplicities}

\author{ Vincenz Busch, Jan Steffen M\"uller\thanks{Supported by DFG-grant STO 299/5-1}}

\newtheorem{thm}{Theorem}
\newtheorem{prop}[thm]{Proposition}
\newtheorem{lemma}[thm]{Lemma}
\theoremstyle{remark}
\newtheorem{rk}[thm]{Remark}
\renewenvironment{proof}{\par\pagebreak[2]\noindent{\it Proof: }}{
 \hfill $\Box$ \medskip}

\newcommand\Q{\mathbb{Q}}
\newcommand\C{\mathbb{C}}

\newcommand\Z{\mathbb{Z}}
\newcommand\R{\mathbb{R}}
\newcommand\Spec{\mathop{\rm Spec}\nolimits}

\renewcommand\O{\mathcal{O}}

\newcommand{\Div}{\operatorname{Div}}
\newcommand{\supp}{\operatorname{supp}}

\newcommand{\SilB}{\operatorname{SilB}}
\newcommand{\To}{\longrightarrow}
\newcommand{\BP}{{\mathbb P}}

\newcommand{\BD}{\textbf{D}}
\newcommand{\OP}{\textbf{P}}
\newcommand{\OO}{\textbf{O}}
\newcommand{\KC}{\mathcal{C}}

\begin{document}
\maketitle

\begin{abstract}
In this short note we prove a formula for local heights on elliptic curves over number fields in terms of intersection theory on a regular model over the ring of integers.
\end{abstract}
\section{Introduction}
Let $K$ be a number field and let $E$ be an elliptic curve in Weierstra\ss\  form defined over $K$. Let $M_K$ denote the set of places on $K$, normalized to satisfy the product formula. For each $v\in M_K$ we denote the completion of $K$ at $v$ by $K_v$ and we let $n_v=[K_v:\Q_v]$ be the local degree at $v$. Then there are certain functions $\lambda_v:E(K_v)\to\R$, called {\em local heights}, such that the canonical height $\hat{h}$ on $E$ can be decomposed as
\begin{equation}\label{CanHtSumLoc}
 \hat{h}(P)=\frac{1}{[K:\Q]}\sum_{v\in M_K}n_v\lambda_v(P).
\end{equation}
In Section \ref{lhts} we discuss our normalization of the local height.

Let $R$ be the ring of integers of $K$ and let $\mathcal{C}$ be the minimal regular model of $E$ over $\Spec(R)$.
If $Q\in E(K)$, we let $\textbf{Q}\in\Div(\mathcal{C})$ denote the closure of $(Q)\in\Div(E)(K)$ and extend this to the group $\Div(E)(K)$ of $K$-rational divisors on $E$ by linearity.

For any non-archimedean $v$ and any divisor $D\in\Div(E)(K)$ of degree zero, Lemma \ref{ExtendD} guarantees the existence of a $v$-vertical $\Q$-divisor $\Phi_v(D)$ on $\mathcal{C}$ such that 
\begin{equation}\label{Phi}
	(\BD+\Phi_v(D)\,.\,F)_v=0\quad\textrm{ for any }v\textrm{-vertical }\Q\textrm{-divisor }F\textrm{ on }\mathcal{C},
\end{equation}
where
$(\cdot \,.\,\cdot)_v$ denotes the intersection multiplicity on $\mathcal{C}$ above $v$. 

In Section \ref{Main-Theorem} we will prove the following result, which is a local analogue of the classical Theorem \ref{fh}.
\begin{thm}\label{main}
 Let $v$ be a non-archimedean place of $K$ and $P\in E(K)\setminus\{O\}$. Suppose that $E$ is given by a Weierstra\ss\ equation that is minimal at $v$. Then we have
 \[
  \lambda_v(P)=2\left(\OP\,.\,\OO\right)_v-\left(\Phi_v((P)-(O)) \,.\,\OP-\OO\right)_v,
 \]
where $\Phi_v((P)-(O))$ is any vertical $\Q$-divisor such that \eqref{Phi} holds for $D=(P)-(O)$.
\end{thm}

Theorem \ref{main} gives a finite closed formula for the local height that is independent of the reduction type of $E$ at $v$. We hope that we can generalize Theorem \ref{main} as described in Section \ref{outlook}.

The first author would like to thank the hospitality of the University of Bayreuth, where most of the research for this paper was done.
\section{Local heights}\label{lhts}

For each non-archimedean place $v$ we let $v:K_v\To \Z\cup\{\infty\}$ denote the surjective discrete valuation corresponding to $v$ and we denote the ring of integers of $K_v$ by $\O_v$.

If $A$ is an abelian variety defined over $K$ and $D$ is an ample symmetric divisor on $A$, one can define the canonical height (or N\'eron-Tate height) $\hat{h}_D$ on $A$ with respect to $D$. In the case of an elliptic curve $A=E$ in Weierstra\ss\ form we use the canonical height $\hat{h}=\hat{h}_{2(O)}$ with respect to the divisor $2(O)$, where $O$ is the origin of $E$.  

For each place $v$ of $K$ there is a local height (or N\'eron function) $\lambda_{D,v}:A(K_v)\to\R$, uniquely defined up to a constant, such that $\hat{h}_D$ can be expressed as a sum of local heights as in \eqref{CanHtSumLoc}, see \cite{LangFund}. For an account of the different normalizations of the local height see \cite[\S4]{cps}; our normalization will correspond to the one used there, so in particular we have
\begin{equation}\label{norm}
\lambda_v(P)=2\lambda^{\SilB}_v(P)+\frac{1}{6}\log|\Delta|_v
\end{equation}
where $\lambda^{\SilB }_v$ is the normalization of the local height with respect to $D=(O)$ used in Silverman's second book \cite[Chapter~VI]{ATAEC} on elliptic curves. 

If $v$ is an archimedean place, then we have a classical characterization of the local height. It suffices to discuss the case $K_v=\C$; here we consider the local height $\lambda':=\lambda_v^{\SilB }$ on $E(\C)\cong\C/\Z \oplus \tau\Z$, where $\text{Im}(\tau) >0$. We set $ q = \exp(2 \pi i \tau)$ and denote by
\[ B_2(T) = T^2-T+\frac{1}{6} \]
	the second Bernoulli polynomial. If $P\in E(\C)\setminus\{O\}$, then we have
 	\[ \lambda'(P) = -\frac{1}{2} B_2 \left(\frac{\text{Im} z}{\text{Im}{\tau}}\right) \log |q| - \log|1-q| - \sum_{n\geq 1} \log |(1-q^nu)(1-q^nu^{-1})| \]
	for any complex uniformisation $z$ of $P$ and $u = \exp(2 \pi i z)$. This is \cite[Theorem~VI.3.4]{ATAEC} and the following result is \cite[Corollary VI.3.3]{ATAEC}:

	\begin{prop}\label{SilvProp}
	For all $P,Q\in E(\C)$ such that $P,Q,P\pm Q\ne O$ we have
	 \[
	  \lambda'(P+Q)+\lambda'(P-Q)=2\lambda'(P)+2\lambda'(Q)-\log|x(P)-x(Q)|+\frac{1}{6}\log|\Delta|.
	 \]
	\end{prop}

	If $v$ is a non-archimedean place we use a Theorem due to N\'eron which concerns the interplay of the local height $\lambda_v$ and the N\'eron model $\mathcal{E}$ of $E$ over $\Spec(\O_v)$. Recall that $\mathcal{E}$ can be obtained by discarding all non-smooth points from $\mathcal{C}\times\Spec(\O_v)$.	Let $(\cdot \,.\,\cdot)_v$ denote the intersection multiplicity on $\mathcal{C}\times\Spec(\O_v)$.

Let $\mathcal{E}_v$ denote the special fiber of $\mathcal{E}$ above $v$; then $\mathcal{E}_v$ has components $\mathcal{E}^{0}_v,\ldots \mathcal{E}^{r}_v$, where $r$ is a nonnegative integer and $\mathcal{E}^{0}_v$ is the connected component of the identity.

	For a prime divisor $D\in \Div(E)(K_v)$ we write its closure in $\mathcal{E}$ as $\BD$ and we extend this operation to $\Div(E)(K_v)$ by linearity. 
	The following proposition is a special case of \cite[Theorem 5.1]{LangFund}:
	\begin{prop}(N\'eron)\label{lang} 
	Let $D\in \Div(E)(K_v)$ and let $\lambda_{D,v}$ be a local height with divisor $D$. For each component $\mathcal{E}_v^{j}$ there is a constant $\gamma_{j,v}(D)$ such that for all
	$ P\in E(K_v)\setminus\supp(D)$
	mapping into $\mathcal{E}^{j}_v$ we have
	\[
	 \lambda_{D,v}(P)=(\BD \,.\,\OP)_v+\gamma_{j,v}(D).
	\] 
	\end{prop}

	\section{Arithmetic intersection theory}\label{Falt-Hril}
	 In this section we briefly recall some basic notions of Arakelov theory on $\mathcal{C}$ and its relation to canonical heights, following essentially \cite{LangArak}. 

		There exists an intersection pairing
		\[ (\cdot\, .\,\cdot ) : \Div(\KC) \times \Div(\KC) \rightarrow \R, \]
	called the {\em Arakelov intersection pairing},
		which, for $D,D'\in\Div(\KC)$ without common component decomposes into
		\[(D\, .\,D' ) = \sum_{v\in M_K}(D\, .\,D' )_v. \]
		In the non-archimedean case $(D\,.\,D')_v$ is the usual intersection multiplicity on $\KC$ above $v$ (defined, for example in \cite[III.\S2]{LangArak}). If $v$ is archimedean, let $g_{D,v}$ denote a {\em Green's function} with respect to $D\times_v\C$ on the Riemann surface $E_v(\C)$ (see \cite[II,\S1]{LangArak}). Then $(D\,.\,D')_v$ is given by $g_{D,v}(D'):=\sum_i n_ig_{D,v}(Q_i)$ if $D'\times_v\C=\sum_i n_iQ_i$. See \cite[IV,\S1]{LangArak}.

Let $v\in M_K$ be non-archimedean.
We say that a divisor $F$ on $\KC$ is {\em $v$-vertical} if $\supp(F)\subset \KC_v$ and we denote the subgroup of such divisors by $\Div_v(\KC)$. We also need to use elements of the group $\Q\otimes\Div_v(\KC)$ of $v$-vertical $\Q$-divisors on $\KC$.

We define the operation $D\to\BD$ on $\Div(E)(K)$ as in Section \ref{lhts}.

		\begin{lemma}(Hriljac)\label{ExtendD}
			For all $D\in\Div(E)(K)$ of degree zero, there exists $\Phi_v(D)\in\Q\otimes\Div_v(\KC)$, unique up to rational multiples of $\KC_v$, such that we have
			\[
				(\BD + \Phi_v(D) \,. \,F)_v=0
			\]
			for any $F\in\Q\otimes\Div_v(\KC)$.
		\end{lemma}
		\begin{proof}
		See for instance \cite[Theorem~III.3.6]{LangArak}.
		\end{proof}
		
Note that we can pick $\Phi_v(D)=0$ if $\mathcal{C}_v$ has only one component.  
This holds for all but finitely many $v$.

	In analogy with a result for elliptic surfaces due to Manin (cf. \cite[Theorem~III.9.3]{ATAEC}), the following theorem relates the Arakelov intersection to the canonical height. See \cite[III,\S5]{LangArak} for a proof.

	\begin{thm}(Faltings, Hriljac) \label{fh}
		Let $D,D' \in \Div(E)(K)$ have degree zero and satisfy $[D]=[D']=P\in \mathrm{Jac}(E)(K)=E(K)$. For each non-archimedean $v$ such that $\mathcal{C}_v$ has more than one component choose some $\Phi_v(D)$ as in Lemma \ref{ExtendD} and set $\Phi(D)=\sum_v \Phi_v(D)$. Then we have
		\[(\BD  + \Phi(D) \,.\,\BD')=  -\hat{h}(P).\]		
	 \end{thm}

	\section{Proof of the Main Theorem}\label{Main-Theorem}
	For a non-archimedean place $v$ we let $E^0(K_v)$ denote the subgroup of points of $E(K_v)$ mapping into the connected component of the identity of the special fiber $\mathcal{E}_v$ of the N\'eron model of $E$ over $\Spec(\O_v)$. We write $\gamma_{j,v}$ for the constant $\gamma_{j,v}(2(O))$ introduced in Proposition \ref{lang} with respect to our local height $\lambda_v$. It is easy to see that our normalization corresponds to the choice $\gamma_{0,v}=0$; therefore we have  
	\begin{equation}\label{lhte0}
	 \lambda_v(P) = (2\OO. \OP)_v =2(\OP. \OO)_v.
	\end{equation}
	 for any $P\in E^0(K_v)\setminus\{O\}$. 
Because $P$ and $O$ reduce to the same component, we also have $\Phi_p((P)-(O))=0$ which proves the theorem for such points.

	Next we want to find the constants $\gamma_{j,v}$ for $j>0$. We will first compare the local height with Arakelov intersections for archimedean places.

	\begin{lemma}\label{GreenLambda1}
	Let $v$ be an archimedean place. The local height $\lambda^{\SilB }_v$ is a Green's function with respect to $D=(O)$ and the canonical volume form on the Riemann surface $E_v(\C)$. Hence the function
	\[
	 g_{P,v}(Q):=\lambda^{\SilB }_v(Q-P)
	\]
	is a Green's function with respect to the divisor $(P)$ for any $P\in E_v(\C)$.
	\end{lemma}
	For a proof see \cite[Theorem II.5.1]{LangArak}. We extend this by linearity to get a Green's function $g_{D,v}$ with respect to any $D\in \Div(E_v)(\C)$.

	\begin{lemma}\label{GreenLambda}
	 Let $v$ be an archimedean place of $K$. For all $P\in E_v(\C)\setminus\{O\}$ and $Q\in E_v(\C)\setminus\{\pm P,O\}$ we have
	 \[
	  g_{D,v}(D_Q)=-\lambda_v(P)-\log|x(P)-x(Q)|_v,
	 \]
	where $D = (P) - (O)$ and $D_Q = ( P + Q ) - (Q)$.
	\end{lemma}
\begin{proof}
We have
\begin{eqnarray*}
 g_{D,v}(D_Q)&=&g_{P+Q,v}(P)-g_{P+Q,v}(O)-g_{Q,v}(P)+g_{Q,v}(O)\\
 &=&2\lambda'(Q)-\lambda'(P+Q)-\lambda'(P-Q),
\end{eqnarray*}
where $\lambda'=\lambda^{\SilB }_v$ and the second equality follows from Lemma \ref{GreenLambda1}. However, by Proposition \ref{SilvProp} we have
\[
 2\lambda'(Q)-\lambda'(P+Q)-\lambda'(P-Q)=-2\lambda'(P)+\log|x(P)-x(Q)|_v-\frac{1}{6}\log|\Delta|_v.
\]
An application of \eqref{norm} finishes the proof of the lemma.
\end{proof}

\begin{lemma}\label{reduction}
	Theorem \ref{main} holds if for each reduction type $\mathcal{K}\notin\{I_0, I_1, II,II^*\}$ there is a prime number $p$ and an elliptic curve $E(\mathcal{K})/\Q$, given by a Weierstra\ss\ equation that is minimal at $p$, satisfying the following conditions:
\begin{itemize}
	\item[(i)] The N\'eron model $\mathcal{E}(\mathcal{K})$ of $E(\mathcal{K})$ has reduction type $\mathcal{K}$ at $p$.
	\item[(ii)] For each connected component $\mathcal{E}(\mathcal{K})_p^j$, there is a point $P_j\in E(\mathcal{K})(\Q)\setminus\{O\}$ reducing to $\mathcal{E}(\mathcal{K})^j_p$.
\end{itemize}
\end{lemma}
\begin{proof}
Let $v$ be a non-archimedean place of $K$, let $k_v$ be the residue class field at $v$. Let $N_v=\frac{n_v}{\log(\#k_v)}$, where $n_v=[K_v:\Q_v]$.	If $P\notin E^0(K_v)$, we have $v_p(x(P))\ge 0$ and hence $(\textbf{P}\,.\, \textbf{O})_p=0$ is immediate.

Now let $\mathcal{K}$ be a reduction type of $E$ at $v$. Then, for any $j\in\{0,\ldots,r\}$, both $\gamma_{j,v}\cdot N_v$ and $\left(\Phi_v((P_j)-(O))\,.\, \OP_j-\OO\right)_v\cdot N_v$ do not depend on $K$, $E$ or $v$, but only on $\mathcal{K}$ and $j$. For the former assertion, see \cite{cps}, where the values of all possible $\gamma_{j,v}$ are determined and for the latter see \cite{cz}. 

Therefore it suffices to show
\begin{equation}\label{jphi}
	\lambda_p(P_j)=\gamma_{j,p}=-\left(\Phi_p((P_j)-(O))\,.\, \OP_j-\OO\right)_p
\end{equation}
for all $j\ne0$, where $P_j\in E(\mathcal{K})(\Q)$ is as in (ii). We can assume $\mathcal{K}\notin\{I_0,I_1,II,II^*\}$, since for those reduction types there is only one connected component.

Let $j\ne0$, let $P=P_j$, let $D = (P) - (O)$ and let $D_Q = ( P + Q ) - (Q)$ for each $Q\in E(\mathcal{K})(\Q)$.

	From Theorem \ref{fh} we deduce
	\[
	 -\sum_{p}(\BD+\Phi_p(D)\,.\, \textbf{D}_\textbf{Q})_p-g_{D,\infty}(D_Q)=\sum_{p}\lambda_p(P)+\lambda_\infty(P)
	\]
	 for any $Q\in E(\mathcal{K})(\Q)\setminus\{P,-P,O\}$. For each prime $p$, the corresponding summand is a rational multiple of $\log p$, so Lemma \ref{GreenLambda} implies 
	\begin{equation}\label{IntLambda}
	 \lambda_p(P)=-(\BD+\Phi_p(D)\,.\, \textbf{D}_\textbf{Q})_p+\log|x(P)-x(Q)|_p
	\end{equation}
	for all primes $p$, by independence of logarithms over $\Q$.

	Now consider $Q=P_0\in E(\mathcal{K})^0(\Q_p)\cap E(\mathcal{K})(\Q)\setminus\{O\}$ and expand
	\[
	 (\BD\,.\, \textbf{D}_\textbf{Q})_p=(\textbf{P}\,.\, \textbf{P+Q})_p-(\textbf{P}\,.\, \textbf{Q})_p-(\textbf{O}\,.\, \textbf{P+Q})_p+(\textbf{O}\,.\, \textbf{Q})_p.
	\]
	By assumption, we have 
	\[
	(\textbf{P}\,.\, \textbf{Q})_p=(\textbf{O}\,.\, \textbf{P+Q})_p=0
	\]
	since the respective points lie on different components. Moreover, because of the N\'eron mapping property, (see \cite[IV,\S5]{ATAEC}) translation by $P$ extends to an automorphism of $\mathcal{E}(\mathcal{K})$, so we have 
	\[
	(\textbf{P}\,.\, \textbf{P+Q})_p=(\textbf{O}\,.\, \textbf{Q})_p=-\frac{1}{2}\max\{v_p(x(Q)),0\}\log p.
	\]
	Since we cannot have $v_p(x(P)-x(Q))>0$, the proof of \eqref{jphi} and hence of the Lemma follows from \eqref{IntLambda}. 
	\end{proof}

	In order to finish the proof of the Theorem, we only need to prove the following result:

	\begin{lemma}
For each reduction type $\mathcal{K}\notin\{I_0, I_1, II, II^*\}$ the elliptic curve $E(\mathcal{K})$ listed in Table \ref{EKs} satisfies the conditions of Lemma \ref{reduction}.
\end{lemma}
\begin{proof}
This is a straightforward check using the proof of Tate's algorithm in \cite[III,\S9]{ATAEC}. If the component group $\Psi(\mathcal{K})$ of $\mathcal{E}(\mathcal{K})$ is cyclic, it suffices to list $P_1\in E(\mathcal{K})(\Q)$ mapping to a generator of $\Psi(\mathcal{K})$ to guarantee the existence of $P_j$ as in Proposition \ref{reduction} for all $j\ne0$. In the remaining case $I^*_n$, $n$ even, we have $\Psi(\mathcal{K})\cong\Z/2\Z\oplus\Z/2\Z$ and hence we need to list two points $P_1$ and $P_2$ mapping to generators of $\Psi(\mathcal{K})$.
\end{proof}

\begin{table}\begin{tabular}{|c|c|c|}
\hline
$\mathcal{K}$ & $p$ & $E(\mathcal{K})$ \\
\hline
\noalign{\hrule height 2pt}
$I_n$ & \multirow{2}{*}{$p>3$} &$y^2=(x+1-p)(x^2-p^{n-1}x+p^n)$\\
\cline{3-3}
$n \ge 2$  & &$P_0=(p-1,0)$;\  $P_1=(p,p)$ \\
\noalign{\hrule height 2pt}
$III$ & \multirow{2}{*}{$7$} & $y^2=x^3+7x+7^2$\\
\cline{3-3}
 	& &$P_0=(-3,1)$;\  $P_1=(0,7)$ \\
\noalign{\hrule height 2pt}
$IV$ & \multirow{2}{*}{$7$} & $y^2=x^3+4\cdot7^2$\\
\cline{3-3}
 	&&$P_0=(-3,13)$;\  $P_1=(0,14)$ \\
\noalign{\hrule height 2pt}
$I_0^*$ & \multirow{2}{*}{$7$} &$y^2 + 7^2 y=x^3 + 7x^2 + 7^2 x$\\
\cline{3-3}
 & &$P_0= (-6,-6);$\  $P_1=(0,0)$;\ $P_2=(14,49)$ \\
\noalign{\hrule height 2pt}
$I_n^*, n\geq 1$ odd          & \multirow{2}{*}{$2$} &$y^2 + 2^k  y = x\cdot (x-(2^k-2))\cdot (x+2^{k+1})$ \\
\cline{3-3}
$n=2k-3$& &$P_0= (-1,2^{k+1}-1);$\  $P_1=(0,0)$ \\
\noalign{\hrule height 2pt}
$I_n^*, n\geq 2$ even & \multirow{2}{*}{$2$} &$y^2 - 2^{k+1}  y = x\cdot (x-(2^k-2))\cdot (x+2^k)$ \\
\cline{3-3}
$n=2k-2$ & &$P_0=(-1,2^k-1)$;\  $P_1=(0,0)$;\ $P_2=(-2^k,0)$\\
\noalign{\hrule height 2pt}
$IV^*$ & \multirow{2}{*}{$7$} & $y^2=x^3+2\cdot7^3x+7^4$\\
\cline{3-3}
 	& &$P_0=(32,239)$;\  $P_1=(0,49)$ \\
\noalign{\hrule height 2pt}
$III^*$ & \multirow{2}{*}{$7$} & $y^2=x^3+7^3x+5\cdot7^5$ \\
\cline{3-3}
 	& &$P_0=(-38,127)$;\  $P_1=(98,1029)$ \\
\noalign{\hrule height 2pt}
\end{tabular}\caption{$E(\mathcal{K})$ for $\mathcal{K}\notin\{I_0, I_1, II, II^*\}$}\label{EKs}\end{table}

\begin{rk}
 It is well-known that $\lambda_v$ is constant on non-identity components of $\mathcal{E}_v$. This follows from Theorem \ref{lang} as above, but we are not aware of any previous result interpreting the constants $\gamma_{j,v}$ in terms of intersection theory. 
\end{rk}
\begin{rk}
It is easy to see that we can consider $P\in E(K_v)$ in the statement of Theorem \ref{main}. In that case, we have to look at the respective Zariski closures on $\KC\times\Spec(\O_v)$ and observe that Lemma \ref{ExtendD} remains correct in the local case.
\end{rk}
\begin{rk} 
Although Theorem \ref{main} requires $E$ to be given by a minimal Weierstra\ss\ equation at $v$, we can find the value of $\lambda_v$ for other models of $E$ using the transformation formula \cite[Lemma~4]{cps}. 
\end{rk}
\begin{rk}
  According to David Holmes, Theorem \ref{main} can also be proved by a direct comparison using N\'eron's original construction of the canonical height pairing. The details will appear in Holmes' forthcoming PhD thesis at the University of Warwick.
\end{rk}

\section{Outlook}\label{outlook}
It would be interesting to generalize Theorem \ref{main} to the case of a Jacobian $J$ of a curve $C$ of genus $g\ge2$. There are analogues of Proposition \ref{lang} in this situation and if we use the divisor $T=\Theta+[-1]^*\Theta$, where $\Theta\in\Div(J)$ is a theta divisor, then Theorem \ref{fh} also generalizes. For instance, if $C$ is hyperelliptic with a unique $K$-rational point $\infty$ at infinity, then every $P\in J(K)$ can be represented using a divisor $D=\sum^d_{i=1}(P_i)-d(\infty)$, where $d\le g$, and a natural analogue of Theorem \ref{main} would be an expression of $\lambda_v=\lambda_{T,v}$ in terms of the intersections $(\OP_i,\bf{\infty})$ and the vertical $\Q$-divisor $\Phi_v(D)$.

This would be interesting, for example, because for $g\ge3$ it is currently impossible to write down non-archimedean local heights explicitly, as one needs to work on an explicit embedding of the Kummer variety $J/\{\pm1\}$ into $\BP^{2^g-1}$ and these become rather complicated as $g$ increases. See \cite[Chapter~4]{jsmthesis} for a discussion. Accordingly, the existing algorithms \cite{David}, \cite{AIntPaper} for the computation of canonical heights use the generalization of Theorem \ref{fh} directly by choosing (rather arbitrarily) divisors $D_1$ and $D_2$ that represent $P$. These algorithms could be simplified significantly if a generalization of Theorem \ref{main} were known.


\end{document}